\documentclass[12pt]{article}

\usepackage{amsmath,amssymb,amsthm}

\let\subset\subseteq \let\eps\varepsilon \let\rho\varrho

\newtheorem{theorem}{Theorem}
\newtheorem{lemma}[theorem] {Lemma}   
\newtheorem{corollary}[theorem] {Corollary}

\newtheorem{definition}[theorem] {Definition} 
\theoremstyle{remark}

\newcommand{\biex}{\textnormal{biex}}
\newcommand{\ex}{\textnormal{ex}}

\title{Dense $H$-free graphs are almost $(\chi(H)-1)$-partite}
\author{Peter Allen\footnote{DIMAP and Mathematics Institute, University of
Warwick, Coventry, CV4 7AL, U.K. Email:~\texttt{allen@ime.usp.br} . Research
supported by the Centre for Discrete Mathematics and its Applications, EPSRC award EP/D063191/1.}}

\begin{document}
\maketitle
\begin{abstract}
By using the Szemer\'edi Regularity Lemma~\cite{Szemeredi76}, Alon and
Sudakov~\cite{AlonSud} recently extended the classical
Andr\'asfai-Erd\H{o}s-S\'os theorem~\cite{AES} to cover general graphs. We prove,
without using the Regularity Lemma, that the following stronger statement is
true.

Given any $(r+1)$-partite graph $H$ whose smallest part has $t$ vertices, there
exists a constant $C$ such that for any given $\eps>0$ and sufficiently large $n$
the following is true. Whenever $G$ is an $n$-vertex graph with minimum degree
\[\delta(G)\geq\left(1-\frac{3}{3r-1}+\eps\right)n~,\] either $G$ contains $H$,
or we can delete $f(n,H)\leq Cn^{2-\frac{1}{t}}$ edges from $G$ to obtain an
$r$-partite graph. Further, we are able to determine the correct order of
magnitude of $f(n,H)$ in terms of the Zarankiewicz extremal function.
\end{abstract}


\section{Introduction}
We define the graph $K_r(s)$ to be the complete $r$-partite graph whose parts
each have $s$ vertices. Given a graph $H$, whose chromatic number is $\chi(H)$,
we examine all the proper $\chi(H)$-colourings of $H$. We choose one whose
smallest colour class is of smallest possible size; then $\sigma(H)$ is the size of
this smallest colour class. Otherwise, our notation is standard.

We recall the classical theorem of Zarankiewicz~\cite{Zaran}:

\begin{theorem} If the $n$-vertex graph $G$ has minimum degree exceeding
$\left(1-\frac{1}{r}\right)n$ then $G$ contains $K_{r+1}$.
\end{theorem}

This theorem is an immediate corollary of Tur\'an's theorem~\cite{Tur}. As is
well known, it is best possible, the extremal example being a complete balanced
$r$-partite graph (sometimes called a Tur\'an graph). An old result of 
Andr\'asfai, Erd\H{o}s and S\'os~\cite{AES}, which amounts to a (very strong)
stability result for Zarankiewicz' theorem, is the following.

\begin{theorem}\label{AESthm} Suppose $r\geq 2$. If the $n$-vertex graph $G$ has
minimum degree exceeding $\left(1-\frac{3}{3r-1}\right)n$ and $G$ does not
contain $K_{r+1}$, then $G$ is $r$-partite.
\end{theorem}

This theorem is best possible; however the extremal example is a little more
complex than the Tur\'an graph. We construct a graph $E_r(n)$ as follows: we
partition $n$ vertices into $r-2$ sets $X_1,\ldots,X_{r-2}$ each containing
$\frac{3n}{3r-1}$ vertices and five sets $Y_1,\ldots,Y_5$ each containing
$\frac{n}{3r-1}$ vertices. Each of these sets is independent; we set every vertex
in each $X_i$ adjacent to all vertices outside $X_i$, and we make
$(Y_i,Y_{i+1\mod 5})$ a complete bipartite graph for each $i$ (so that the five
sets form a blow-up of $C_5$). It is straightforward to check that each vertex
has degree $\left(1-\frac{3}{3r-1}\right)n$; since $\chi(C_5)=3$ the chromatic
number of $E_r(n)$ is $r+1$, but $E_r(n)$ does not contain $K_{r+1}$.

Erd\H{o}s and Stone~\cite{ErdSto} extended Zarankiewicz' theorem, showing that
for any fixed graph $H$, the chromatic number of $H$ governs the minimum degree
threshold at which $H$ appears in a large graph $G$:

\begin{theorem}\label{ErdosStone} Let $H$ be any fixed graph with chromatic
number $r+1$. If the $n$-vertex graph $G$ has minimum degree exceeding
$\left(1-\frac{1}{r}+o(1)\right)n$ then $G$ contains $H$.
\end{theorem}

Although the extremal graphs for this theorem are not necessarily $r$-partite, it
is true that one may delete $o(n^2)$ edges from any extremal graph to obtain an
$r$-partite graph. Indeed, it is not hard to show that there exists
$\rho=\rho(H)>0$ such that deletion of only $O(n^{2-\rho})$ edges from an
extremal graph yields an $r$-partite graph.

Quite recently, Alon and Sudakov~\cite{AlonSud} gave an extension of Andr\'asfai,
Erd\H{o}s and S\'os' result to cover all fixed graphs $H$ (Erd\H{o}s and
Simonovits~\cite{ErdSimo} had previously considered the case when $H$ is
critical, i.e. when there is an edge of $H$ whose removal decreases the chromatic
number):

\begin{theorem}\label{AlS} Let any fixed graph $H$ with chromatic number $r+1$
and constant $\eps>0$ be given. Then there exist $\rho=\rho(H)>0$ and
$n_0=n_0(H,\eps)$ such that the following holds. If $n\geq n_0$ and $G$ is an
$n$-vertex graph with minimum degree exceeding
$\left(1-\frac{3}{3r-1}+\eps\right)n$ which does not contain $H$, then one can
delete at most $O(n^{2-\rho})$ edges from $G$ to yield an $r$-partite graph.
\end{theorem}

Alon and Sudakov gave a value for the constant $\rho(H)$. They showed that if
we have $H\subset K_{r+1}(s)$ then we may take $\rho(H)=1/4r^{2/3}s$. The
purpose of this paper is to give a simpler proof (avoiding the use of the Regularity Lemma) which gives
the correct order of magnitude of the number of edges that must be deleted
(albeit in terms of the Zarankiewicz problem).

Recall that given a family $\mathcal{H}$ of graphs, $\ex(n,\mathcal{H})$ is
defined to be the maximum number of edges in an $n$-vertex graph which does not
contain a copy of any graph $H\in\mathcal{H}$.

Given a graph $H$, we define a quantity $\biex(n,H)$ as follows. Let
$c:V(H)\rightarrow [\chi(H)]$ be any proper $\chi(H)$-colouring of $H$. Let
$S_c=c^{-1}(\{1,2\})$ be the vertices receiving colours $1$ and $2$ in this
colouring. Consider the family of graphs $\mathcal{F}$ containing all graphs
of the form $H[S_c]$ for some proper $\chi(H)$-colouring $c$ of $H$. Then we
set $\biex(n,H)=\ex(n,\mathcal{F})$.

We note that if $H$ is a complete $r$-partite graph, whose smallest part has $t$
vertices and whose next smallest part has $s$ vertices, then
$\biex(n,H)=\ex(n,K_{t,s})$.

The problem of estimating $\ex(n,H)$ when $H$ is bipartite (or, more generally,
$\ex(n,\mathcal{H})$ for a family $\mathcal{H}$ of bipartite graphs) is the
Zarankiewicz problem; for most $H$ it is quite far from being solved. However an upper bound is provided by the following classical theorem of K\"ov\'ari, S\'os and
Tur\'an~\cite{KST}.

\begin{theorem}\label{KSTthm} Let $1\leq t\leq s$ be fixed integers. If $G$ is
any $n$-vertex graph with $\Omega(n^{2-\frac{1}{t}})$ edges, then $G$ contains
$K_{t,s}$.
\end{theorem}

We note that for $t=1,2,3$ and when $s\geq t!+1$ there exist lower bound
constructions matching the upper bound of Theorem~\ref{KSTthm}
(see~\cite{Reiman,Brown,KRS}); for $t\geq 4$ the best known lower bound is
$\Omega(n^{2-\frac{2}{t+1}})$, but it is conjectured that the correct bound is
$\Theta(n^{2-\frac{1}{t}})$.

We can now state our main theorem.

\begin{theorem}\label{AESext} To any graph $H$ with chromatic number $r+1$ there
is associated a constant $C=C(H)$ such that whenever $\eps>0$ is given, there is
$n_0$ for which the following holds. Whenever $n\geq n_0$ and $G$ is an
$n$-vertex graph with minimum degree exceeding
$\left(1-\frac{3}{3r-1}+\eps\right)n$ which does not contain $H$, then one can
delete at most $C\biex(n,H)$ edges from $G$ to obtain an $r$-partite graph.
\end{theorem}

This theorem is best possible up to the value of $C$. For comparison with the
result of Alon and Sudakov, suppose $H\subset K_{t,s,s,\ldots,s}$ has chromatic
number $r+1$, where $t\leq s$. Then, applying Theorem~\ref{KSTthm}, we have
\[\biex(n,H)\leq\ex(n,K_{t,s})=O(n^{2-\frac{1}{t}})~.\]

It follows that if $G$ satisfies the conditions of Theorem~\ref{AESext}, then
Theorem~\ref{AlS} guarantees that $G$ can be made $r$-partite by deleting
$O(n^{2-\frac{1}{4r^{2/3}s}})$ edges; Theorem~\ref{AESext} strengthens this to
$Cn^{2-\frac{1}{t}}$ edges. On the assumption that the conjectured bound in
the Zarankiewicz problem is correct, this is best possible up to the value of
the multiplicative constant. Furthermore, the constant hidden behind the
$O(\cdot)$ notation in Theorem~\ref{AlS} depends upon $\eps$;
specifically, it grows as a polynomial function of $1/\eps$, whereas the
constant $C$ in Theorem~\ref{AESext}, while surely much larger than it
`should' be, does not depend on $\eps$. Finally, owing to the use of the
Regularity Lemma, the constant $n_0$ in Theorem~\ref{AlS} has an exceptionally
unpleasant dependence on $\eps$, $r$ and $s$.

We give two constructions which demonstrate the tightness of our theorem.

Given $H$, let $E$ be an $n$-vertex graph with $\biex(n,H)$ edges and not
containing any of the forbidden bipartite subgraphs. Let $E'$ be an $n/r$-vertex
bipartite subgraph of $E$ containing the maximum possible number of edges. Note
that $e(E')> e(E)/4r^2=\Omega(\biex(n,H))$.

Consider the graph $G$ obtained from the complete balanced $r$-partite graph by
replacing one part with $E'$. This graph has minimum degree $\frac{r-1}{r}n$, and
does not contain a copy of $H$. However to make $G$ $r$-partite we must delete 
$\Omega(\biex(n,H))$ edges.

Alon and Sudakov asked whether it is possible to replace the term $\eps n$ in the
minimum degree of their theorem with an $O(1)$ term. It is not possible; indeed,
for any $\mu>0$ there are graphs $H$ such that the corresponding term must be
larger than $n^{1-\mu}$.

Consider the following modification of $E_r(n)$. Let $c$ be some sufficiently
small positive quantity.  We let each of
the independent sets $Y_1,\ldots,Y_5$ have $\frac{n}{3r-1}+(r-2)cn^{1-2/t}$
vertices. We let each of the independent sets $X_1,\ldots,X_{r-2}$ have
$\frac{3}{3r-1}-5cn^{1-2/t}$ vertices. Finally, we take a $K_{t,t}$-free graph $E'$
on $|Y_1|$ vertices with minimum degree $(3r-1)cn^{1-2/t}$: provided $c>0$ is
chosen sufficiently small, such a graph exists. We replace each of the
independent sets $Y_1,\ldots,Y_5$ with $E'$ to obtain $E'_{r,t}(n)$. Now
observe that the minimum degree of $E'_{r,t}(n)$ is
$\frac{3r-4}{3r-1}n+5cn^{1-2/t}$. However it is not possible to find a copy of
$K_{r+1}(2t)$ in $E'_{r,t}(n)$. The reason is that it would be necessary to
find a copy of $K_3(2t)$ within the graph induced by $Y_1\cup\ldots\cup
Y_5$; this would require that one of the $Y_i$ contained $K_{t,t}$, which
by construction is false. Finally,
it is clear that to make $E'_{r,t}(n)$ $r$-partite requires the removal of
$\Omega(n^2)$ edges.

\section{Constructing $(r+1)$-partite graphs}

Given an $(r+1)$-partite graph $H$, a large graph $G$, and a family $\mathcal{F}$
consisting of the bipartite subgraphs of $H$ whose removal decreases the
chromatic number of $H$ by two, we describe a construction of the graph $H$ from
a suitably well-structured set of copies of $K_{r+1}$ in $G$. Alon and Sudakov
made use of a related construction: the difference is that their construction as
its first step finds (by use of the K\"ov\'ari-S\'os-Tur\'an theorem) one
specific bipartite subgraph of $G$ and proceeds to build $H$ using it. Our
construction avoids this, relying instead on counting the number of suitable
objects until the final step in the construction. This difference is primarily
responsible for our improved bounds.

Given a graph $G$ and a vertex $v\in G$, let $G_v$ be the \textit{neighbourhood
graph} obtained by deleting from $G$ every edge which is not contained in the
neighbourhood of $v$.

We give first a counting variant of a lemma of Erd\H{o}s~\cite{Erd1964}; this
is essentially a statement about dense hypergraphs generalising the
K\"ov\'ari-S\'os-Tur\'an theorem.

\begin{lemma}~\label{ErdCount} For every $r$, $s$ and $\eps>0$ there exists
$\delta=\delta_{r,s}(\eps)>0$ such that the following holds for sufficiently
large $n$. If the $n$-vertex graph $G$ contains at least $\eps n^r$ copies of
$K_r$, then $G$ contains $\delta_{r,s}(\eps) n^{rs}$ copies of $K_r(s)$.
\end{lemma}
\begin{proof}
For $r=1$ the statement holds trivially. We complete the proof by induction.

Let $G$ be an $n$-vertex graph containing $\eps n^r$ copies of $K_r$: then there
are some $\eps n/2$ vertices $D$ of $G$ which are each contained in $\eps
n^{r-1}/2$ copies of $K_r$ in $G$. By construction, for each $d\in D$, $G_d$
contains $\eps n^{r-1}/2$ copies of $K_{r-1}$; by induction it contains
$\delta_{r-1,s}(\eps/2)n^{(r-1)s}$ copies of $K_{r-1}(s)$.

For a given copy $S$ of $K_{r-1}(s)$, let $d_S$ be the number of vertices of $D$
whose neighbourhoods contain $S$. Then we have (using the convention that
$\binom{a}{b}=0$ when $a<b$) at least $\frac{1}{r}\sum_S \binom{d_S}{s}$ copies
of $K_r(s)$ contained in $G$. Since the mean value of $d_S$ is at least
$\delta_{r-1,s}(\eps/2)|D|$, applying Jensen's inequality the number of copies of
$K_r(s)$ in $G$ is at least \[\frac{1}{r}\sum_S\binom{d_S}{s}\geq\frac{1}{r}
\delta_{r-1,s}(\eps/2)n^{(r-1)s}\binom{\delta_{r-1,s}(\eps/2)|D|}{s}=\delta_{r,s}(\eps)n^{rs}~,\]
as required.
\end{proof}

Note that the value of $\delta_{r,s}(\eps)$ obtained by the above method is
polynomial in $\eps$. 

To complete our construction, we give the following corollary of Lemma~\ref{ErdCount}.

\begin{corollary}\label{NiceCons} Given $\eps>0$ and $H$ there exists $C$ such
that for sufficiently large $n$ the following is true. Every $n$-vertex graph
$G$ in which there are more than $C\biex(n,H)$ edges $E$ of $G$, each contained in
$\eps n^{r-1}$ copies of $K_{r+1}$, contains $H$.
\end{corollary}
\begin{proof} Let $G$ be a graph with a set $E$ of edges each of whose common neighbourhoods contains
$\eps n^{r-1}$ copies of $K_{r-1}$. Suppose that $n$ is large enough to permit
us to conclude, by Lemma~\ref{ErdCount}, that the common neighbourhood of each
edge of $E$ contains at least 
$\delta_{r-1,v(H)}(\eps) n^{(r-1)v(H)}$ copies of $K_{r-1}(v(H))$. Let
$C=1/\delta_{r-1,v(H)}(\eps)$. Suppose furthermore that $|E|>
C\biex(n,H)$.

 By averaging, there is one copy $S$ of $K_{r-1}(v(H))$ in $G$
which lies in the common neighbourhood of each of the edges $E'\subset E$,
with $|E'|>\biex(n,H)$. By definition of $\biex(n,H)$, the edges $E'$ must
contain a copy of some bipartite subgraph of $H$ in $\mathcal{F}$. Let this
subgraph be $B$. Then $B\cup S$ contains $H$. 
\end{proof}

Note that the value of $\delta_{r,s}(\eps)$ given by Lemma~\ref{ErdCount} is
clearly far smaller than the truth; but this affects only the constant $C$;
furthermore, the dependence on $\eps$ is polynomial.

\section{Proof of Theorem~\ref{AESext}}

We first prove a density version of Theorem~\ref{AESthm}. We note that Alon and
Sudakov~\cite{AlonSud} proved a similar lemma; however their method (while in
most ways similar to ours) obtained a first `coarse' version by application of
the Szemer\'edi Regularity Lemma. We avoid this by making use of an induction
argument.

\begin{lemma}\label{AESDens} Given $r$ and $\eps$, let $\mu=\eps^r/r!$ and
$\eta=\eps^{r+1}/(r+1)!$. Then whenever $n$ is sufficiently large, the following is true. Any
$n$-vertex graph $G$ with $\delta(G)>\left(1-\frac{3}{3r-1}+4\eps\right)n$
either contains more than $\eta n^{r+1}$ copies of $K_{r+1}$, or has a
partition into $D\cup V_1\cup\ldots\cup V_r$, with the properties that
$\Delta(G[V_i])\leq\eps n$ for each $i$, each vertex of $D$ is contained in
more than $\mu n^r$ copies of $K_{r+1}$, and $|D|\leq \eps n$.
\end{lemma}

Note that when $\eps=0$ we have $\mu=\eta=0$, and we obtain the statement
of Theorem~\ref{AESthm}. The intuition is that since we are looking at graphs
which do not contain a high density of copies of $K_{r+1}$, rather than not
containing any at all, we must expect that there may be some small set of
vertices, and a few edges leaving every vertex, which `misbehave'. These are,
respectively, the set $D$ and the replacement of the independent sets of
Theorem~\ref{AESthm} with sets which simply have restricted maximum degree.

\begin{proof} We prove the lemma by induction. The $r=1$ case is a triviality: either there are
more than $\eps n$ vertices of degree exceeding $\mu n$, in which case $G$
certainly contains more than $\eta n^2$ edges, or we can let $D$ be the set
of all vertices of degree exceeding $\mu n$, and together with
$V_1=V(G)\setminus D$ the partition conclusion is satisfied.

Suppose $r\geq 2$. We assume as our induction hypothesis that the lemma holds
for $r-1$.

Let $G$ be an $n$-vertex graph with minimum degree
$\left(1-\frac{3}{3r-1}+4\eps\right)n$. We presume $G$ contains at most $\eta
n^{r+1}$ copies of $K_{r+1}$.

Let $D\subset V(G)$ be the set of all
vertices $d\in G$ such that there are more than $\mu n^r$ copies of
$K_r$ in $\Gamma(d)$. Then $|D|\leq\eps n$ since $G$ contains at most $\eta
n^{r+1}$ copies of $K_{r+1}$.

Let $G'=G[V(G)-D]$. This graph has minimum degree greater than
$\left(\frac{3r-4}{3r-1}+3\eps\right)n$; none of its vertices are
contained in more than $\mu n^r$ copies of $K_r$.

Let $X_1$ be a maximum cardinality set in $V(G')$ with the property that
$\Delta(G'[X_1])\leq\eps n$. Let $v\in X_1$.

Consider the graph $N=G'[\Gamma(v)\setminus X_1]$. Because $v\notin D$, the
neighbourhood graph $G_v$ contains at most $\mu n^r$ copies of $K_r$, and so in
particular $N$ contains at most $\mu n^r$ copies of $K_r$. Because
$\Delta(G'[X_1])\leq \eps n$, $v(N)>\frac{3r-4}{3r-1}n+2\eps n$. Now consider
$u\in N$. We have
\[d_N(u)>v(N)-\left(\frac{3}{3r-1}-4\eps\right)n\]
\[>v(N)-\left(\frac{3}{3r-1}-4\eps\right)\frac{3r-1}{3r-4}v(N)>\left(\frac{3r-7}{3r-4}+4\eps\right)v(N)~.\]

By induction, we have that $N$ has a partition $V(N)=B\cup X_2\cup\ldots\cup
X_r$, where $|B|\leq\eps n$ and $\Delta(N[X_i])\leq\eps n$ for each of the $r-1$
sets $X_2,\ldots,X_r$.

Because $X_1$ has maximum cardinality subject to $\Delta(G'[X_1])\leq\eps n$,
$|X_1|\geq |X_i|$ for each $i$. In particular, we have
\[|X_1|+\ldots+|X_r|\geq \left(\frac{3r-4}{3r-1}+\eps
\right)\frac{rn}{r-1}\geq \frac{(3r-4)rn}{(3r-1)(r-1)}+\eps n~.\]

Since every vertex in $G$ has more than $\frac{3r-4}{3r-1}n+4\eps n$ neighbours
in $G$, and since for each $i$ we have $\Delta(G[X_i])\leq\eps n$, it follows
that $|X_i|<\frac{3}{3r-1}n$ for each $i$.

Now suppose that for some $i$ we have $|X_i|\leq\frac{2}{3r-1}n$. Because $X_1$
was chosen to be maximal, we may assume $2\leq i\leq r$; without loss of
generality let us suppose $i=r$. We have
$|B|+|X_2|+\ldots+|X_r|=v(N)\geq\frac{3r-4}{3r-1}n+2\eps n$, and since also
$|B|\leq\eps n$, we have $|X_2|+\ldots+|X_{r-1}|\geq\frac{3r-6}{3r-1}n+\eps n$.
It follows that among the $r-2$ sets $X_2,\ldots,X_{r-1}$, there must be one whose size exceeds
$\frac{3r-6}{(3r-1)(r-2)}n=\frac{3}{3r-1}n$, which is a contradiction. Thus we
have that for each $i$, $\frac{2}{3r-1}n<|X_i|<\frac{3}{3r-1}n$.

Now, if we have any two adjacent vertices $u$ and $v$ of $G'$ whose codegree
exceeds $\frac{3r-6}{3r-1}n+\eps n$, then we may construct a clique $K_{r+1}$
extending $uv$ greedily by simply picking any common neighbour of the so far
chosen vertices at each step. At the final step (and therefore at all steps)
we have at least $\eps n$ choices. It follows that any edge $uv$ of $G$ in
which the common neighbourhood of $u$ and $v$ exceeds $\frac{3r-6}{3r-1}n+\eps
n$ lies in more than $\eps^{r-1}n^{r-1}/(r-1)!$ cliques $K_{r+1}$.

Furthermore, if $u$ has more than $\eps n$ neighbours with each of which its
codegree exceeds $\frac{3r-6}{3r-1}n+\eps n$, then $u$ lies in more than
$\eps^r n^r/r!=\mu n^r$ copies of $K_{r+1}$. This contradicts $u\notin D$.

Since $\Delta(G[X_i])\leq \eps n$, if a vertex $u$ outside $X_i$ has less than
$|X_i|-\frac{n}{3r-1}$ neighbours in $X_i$, then the codegree of $u$ and any
neighbour $v\in X_i$ exceeds $\frac{3r-6}{3r-1}n+\eps n$. It follows that any
vertex of $G'$ outside $X_i$ has either fewer than $\eps n$ neighbours in $X_i$
or more than $|X_i|-\frac{n}{3r-1}$ neighbours in $X_i$.

Consider the set $L_i$ of vertices of $L$ which all have less than $\eps n$
neighbours in $X_i$. Any one of these vertices has codegree exceeding
$\frac{3r-6}{3r-1}n+\eps n$ with any other, and with any vertex of $X_i$. It
follows that $L_i\cup X_i$ has maximum degree $\eps n$. Let this set be
$V_i$. Let the vertices of $G'$ not in any $X'_i$ be $L'$.

If $L'=\emptyset$ then we have $V(G)=D\cup V_1\cup\ldots\cup V_r$ is the
desired partition. So we may assume there is a vertex $l\in L'$. This vertex
is non-adjacent to fewer than $\frac{n}{3r-1}$ vertices of each set
$V_i$. It is convenient to assume that the sets $V_1,\ldots,V_r$ are in order of
decreasing size.

Finally, consider the following greedy construction. We start with the vertex
$l\in L'$. We now choose vertices $v_1,\ldots,v_r$ from the
respective sets $V_1,\ldots,V_r$, such that after each choice the vertices
chosen together with $l$ form a clique.

At the first step we have more than $|V_1|-\frac{n}{3r-1}$ choices for $v_1$. At
the second step we have more than
\[|V_2|-\frac{n}{3r-1}-\left(\frac{3}{3r-1}-4\eps\right)n+(|V_1|-\eps
n)=|V_1|+|V_2|-\frac{4}{3r-1}n+3\eps n\] choices for $v_2$; there are less than
$\frac{n}{3r-1}$ non-neighbours of $l$ in $V_2$, and at most
$\frac{3n}{3r-1}-4\eps n$ non-neighbours of $v_1$ in $G$, of which at least
$|V_1|-\eps n$ are in $V_1$. In general, for each $2\leq i\leq r$, we have at the
$i$th step more than \[|V_1|+\ldots+|V_i|-\frac{3i-2}{3r-1}n+3\eps n\] choices
for $v_i$. Because the sets $V_1,\ldots,V_r$ are in order of decreasing size, the
number of choices is least when choosing either $v_1$ or $v_r$. Since $|V_1|\geq
|X_1|>\frac{3r-4}{(3r-1)(r-1)}n\geq\frac{2}{3r-1}n$, the number of choices for
$v_1$ is greater than $\frac{n}{3r-1}$. Since
\[|V_1|+\ldots+|V_r|\geq
|X_1|+\ldots+|X_r|\geq\frac{(3r-4)r}{(3r-1)(r-1)}n+\eps n~,\] the number of
choices for $v_r$ is at least $\frac{r-2}{(3r-1)(r-1)}n+4\eps n$. It follows that at each
step there are more than $\eps n$ choices; therefore $l$ is contained in more
than $\eps^rn^r\geq\mu n^r$ copies of $K_{r+1}$ in $G$, which contradicts
$l\notin D$.
\end{proof}

At last, we can complete the proof of our main theorem. Again, our method is
similar to that of Alon and Sudakov~\cite{AlonSud}; we take a little more care
in order to ensure that the constant $C$ in our theorem is independent of
$\eps$.

\begin{proof}[Proof of Theorem~\ref{AESext}]
Given $r\geq 2$ and $\eps>0$, let $G$ be a sufficiently large $n$-vertex graph
with $\delta(G)\geq\left(1-\frac{3}{3r-1}+\eps\right)n$ which does not contain the
$(r+1)$-partite graph $H$.

By Lemma~\ref{AESDens} there exist positive constants
$\eta,\mu$ such that either $G$ contains $\eta n^r$ copies of $K_{r+1}$ or
$V(G)$ may be partitioned as $V(G)=D\cup V_1\cup\ldots\cup V_r$, where
$\Delta(G[V_i])\leq\eps n/4$ for each $i$, each vertex of $D$ is contained
in at least $\mu n^r$ copies of $K_{r+1}$, and $|D|\leq\eps n/4$.

When $n$ is sufficiently large, by Lemma~\ref{ErdCount} every graph $G$ with
$\eta n^{r+1}$ copies of $K_{r+1}$ contains $K_{r+1}(v(H))$ and thus $H$. It
follows that $V(G)$ possesses the given partition.

As in the proof of Lemma~\ref{AESDens}, for each $i$, since $\Delta(V_i)\leq\eps
n/4$ and $\delta(G)>\frac{3r-4}{3r-1}n+\eps n$, we have
$|V_i|<\frac{3}{3r-1}n-3\eps n/4$. Again, if for some $i$ we have
$|V_i|\leq\frac{2}{3r-1}n$ then among the $r-1$ sets $V_1,\ldots,V_r$ remaining
there must be one whose size is at least \[\left(n-\eps
n/4-\frac{2}{3r-1}n\right)/(r-1)>\frac{3}{3r-1}n-\eps n/2~,\] which again is a contradiction. Thus for each $i$ we have
$\frac{2}{3r-1}n<|V_i|<\frac{3}{3r-1}n$.

We alter slightly the partition given by Lemma~\ref{AESDens} as follows. For
each $1\leq i\leq r$, let $W_i$ be the set of vertices with at most $\frac{n}{4(3r-1)}$
neighbours in $V_i$. Let $Y_i$ be the vertices of $D$ with more than $\frac{n}{4(3r-1)}$
neighbours, but less than $|V_i|-\frac{3}{2(3r-1)}n$ neighbours in $V_i$. Let $X$ be the vertices of
$D$ not contained in any set $W_i$ or $Y_i$. By definition of $V_i$, we have
$V_i\subset W_i$ for each $i$.

Consider the vertex $x\in X$. We make use of a greedy construction as in the proof
of Lemma~\ref{AESDens}. We presume that the sets $V_1,\ldots,V_r$ are in order
of decreasing size. We choose greedily vertices $v_1,\ldots,v_r$ in sets
$V_1,\ldots,V_r$ (in that order), such that the set $\{x,v_1,\ldots,v_r\}$ are
the vertices of an $(r+1)$-clique in $G$. As in the proof of
Lemma~\ref{AESDens}, at the $i$th step we have at least
\[|V_1|+\ldots+|V_i|-\frac{3}{2(3r-1)}n-\frac{3i-3}{3r-1}n+3\eps n/4\]
choices for $v_i$. As before, since the sets $V_i$ are in order of decreasing
size the number of choices is fewest at either the first or the last step. The
number of choices at the first step is at least
$|V_1|-\frac{3}{2(3r-1)}>\frac{1}{2(3r-1)}n$; since the sets
$V_1,\ldots,V_r$ together cover all of $G$ except the at most $\eps n/4$
vertices of $D$, the number of choices at the last step is at least
\[n-\eps n/4-\frac{3}{2(3r-1)}n-\frac{3r-3}{3r-1}n+3\eps
n/4>\frac{1}{2(3r-1)}n~.\]

It follows that at every step there are at least $\frac{1}{2(3r-1)}n$
choices, and hence $x$ is contained in at least
\[\left(\frac{n}{2(3r-1)}\right)^r\]
copies of $K_{r+1}$ in $G$.

Consider the vertex $y\in Y_i$. Let $u$ be any neighbour of $y$ in
$V_i$. The common neighbourhood of $u$ and $y$ contains at least
\[2\left(\frac{3r-4}{3r-1}+\eps \right)n-\left(n-\frac{3}{2(3r-1)}n+\eps
n/4\right)>\frac{6r-11}{2(3r-1)}n\]
vertices. Now we construct an $(r+1)$-clique greedily starting from $uy$. At
the final step, and thus at every step, we have at least $\frac{n}{2(3r-1)}$
choices. It follows that $uy$ lies in at least
$\left(\frac{n}{2(3r-1)}\right)^r/(r-1)!$ copies of $K_{r+1}$ in $G$. Since
$y$ has at least $\frac{n}{4(3r-1)}$ neighbours in $V_i$, $y$ lies in at least
$\left(\frac{n}{4(3r-1)}\right)^r/r!=\gamma n^r$ copies of $K_{r+1}$ in $G$.

Finally we have that every vertex of $Z=Y_1\cup\ldots\cup Y_r\cup X$ lies in at
least $\gamma n^r$ copies of $K_{r+1}$ in $G$.

Now by Lemma~\ref{ErdCount} there exists
$\delta>0$ such that whenever $n$ is sufficiently large, every graph $G$
with $\gamma n^r$ copies of $K_r$ contains $\delta n^{rv(H)}$ copies of
$K_r(v(H))$. If $|Z|>(\sigma(H)-1)/\delta$, then there is one copy $S$ of
$K_r(v(H))$ in $G$ which is in the neighbourhood of each of $\sigma(H)$
vertices $B$ of $G$. But then $H\subset G[B\cup S]$, which is a
contradiction. It follows that $|Z|\leq(\sigma(H)-1)/\delta$. It is important
to note that $\gamma$, and hence $\delta$, are independent of $\eps$.

Finally, let $E$ be the set of edges of $G$ which are contained in any one of
the sets $W_i$.

For any edge $uv\in E$, there is $i$ such that $u,v\in W_i$. Then the common
neighbourhood of $u$ and $v$ in $V(G)$ contains at least
\[2\left(\frac{3r-4}{3r-1}+\eps\right)n-\left(n-|V_i|+\frac{n}{2(3r-1)}\right)\geq
\frac{6r-11}{2(3r-1)}n+2\eps n\] vertices, since both $u$ and $v$ are adjacent
to at most $\frac{n}{4(3r-1)}$ vertices of $V_i$. As before, we can extend $uv$ to a
clique $K_{r+1}$ by choosing vertices greedily; at each stage we have at least
$\frac{n}{2(3r-1)}$ choices, and hence $uv$ is contained in at least
$\frac{n^{r-1}}{(6r-2)^{r-1}(r-1)!}$ copies of $K_{r+1}$. By
Corollary~\ref{NiceCons}, since $G$ does not contain $H$, there exists $C'$ such
that $|E|\leq C'\biex(n,H)$. Observe that $C'$ does not depend on $\eps$.

If $\biex(n,H)<n-1$, then it must be the case that there is some
bipartite subgraph $F$ of $H$ such that $F\subset K_{1,n-1}$ and the graph
$H[V(H)\setminus V(F)]$ is $(r-1)$-colourable. But then there is a proper
$(r+1)$-colouring of $H$ in which one colour class has size one; so $\sigma(H)=1$.

Upon deleting from $G$ all edges incident to $Z$ or contained in $E$, one
obtains an $r$-partite graph. The total number of edges deleted is at most
$n(\sigma(H)-1)/\delta+C'\biex(n,H)$. Since $n|Z|>0$ only if $\sigma(H)>1$, i.e. only if
$\biex(n,H)\geq n-1$, we have $n|Z|+C'\biex(n,H)\leq C\biex(n,H)$,
and $C$ is as required independent of $\eps$ since $C'$ and $\delta$ are. 
\end{proof}

\section{Concluding remarks}

Perhaps the main conclusion of this paper is that (if such is necessary) there
is a further motivation for solving the Zarankiewicz problem of determining
$\ex(n,\mathcal{F})$ for all families $\mathcal{F}$ of bipartite graphs.

However there remain some open questions which are independent of the
Zarankiewicz problem.

First, it would be interesting to know what the best possible value of
$\mu(H)$ is such that the following statement is true.

Given $H$, with $\chi(H)=r+1$, there exists $C$ such that for all sufficiently
large $n$, if $G$ is an $n$-vertex $H$-free graph with minimum degree at least
$\frac{3r-4}{3r-1}n+\Theta(n^{1-\mu})$, then $G$ can be made $r$-partite by
deleting at most $C\biex(n,H)$ edges.

It follows (by careful analysis of the proof given) that $\mu(H)$ must always be
positive: but it seems likely that the value so obtained is much smaller than optimal.

Second, although we have shown that the correct
number of edges which we should delete from a dense $H$-free graph $G$ to obtain a
$(\chi(H)-1)$-partite graph is $\Theta(\biex(n,H))$, it seems certain that the
multiplicative constants proved for our upper and lower bounds are not best
possible. We have made no particular effort to optimise our upper bound: but
probably such effort using our techniques would produce only a somewhat less bad upper
bound.

It would be interesting to know whether there exists a best possible
value for the constant $C$, and if so, what it is. It seems likely that
(despite the result of this paper) the best possible value will depend upon
$\eps$.

\section*{Acknowledgement}

The author would like to thank Daniela K\"uhn and Deryk Osthus for suggesting this nice problem.

\providecommand{\bysame}{\leavevmode\hbox to3em{\hrulefill}\thinspace}
\providecommand{\MR}{\relax\ifhmode\unskip\space\fi MR }
\providecommand{\MRhref}[2]{%
  \href{http://www.ams.org/mathscinet-getitem?mr=#1}{#2}
}
\providecommand{\href}[2]{#2}

\appendix
\section{General monotone properties}

The purpose of this appendix is to explain to what extent we can generalise the
preceding results to the setting of considering graphs from any monotone
property $\mathcal{P}$, with any number (finite or infinite) of forbidden subgraphs, as
opposed to the property defined by fixing one graph $H$, and excluding it as a
subgraph. That this should be possible was suggested by Rob Morris; the
mathematical contents of this section arose from discussions with Julia
B\"ottcher, Simon Griffiths and Rob Morris.

Let $\mathcal{P}$ be a monotone property (that is, for every $G\in\mathcal{P}$
and every subgraph $G'$ of $G$, we have $G'\in\mathcal{P}$). We call a graph $L$
a \emph{forbidden graph} for $\mathcal{P}$ if we have $L\notin\mathcal{P}$, and
furthermore for every proper subgraph $L'$ of $L$, $L'\in\mathcal{P}$. Then
there is a family $\mathcal{L}(\mathcal{P})$ consisting of all forbidden
subgraphs of $\mathcal{P}$, and $\mathcal{P}$ consists precisely of all graphs
with no subgraph contained in $\mathcal{L}(\mathcal{P})$; we say $\mathcal{P}$
is the $\mathcal{L}(\mathcal{P})$-free property.

In the language of monotone properties, Theorem~\ref{ErdosStone} is a statement
about the largest possible minimum degree of any $n$-vertex graph in the
$\{H\}$-free property. The corresponding result for general properties was
proved by Erd\H{o}s~\cite{ErdProp} and Simonovits~\cite{SimProp}. Given a set
$\mathcal{L}$ of graphs, let $p=p(\mathcal{L})=\min_{L\in\mathcal{L}}\chi(L)-1$.

\begin{theorem} Let $\mathcal{L}$ be any fixed set of non-empty graphs. If the
$n$-vertex graph $G$ has minimum degree exceeding
$\left(1-\frac{1}{p(\mathcal{L})}+o(1)\right)n$ then $G$ is not
$\mathcal{L}$-free.
\end{theorem}

A natural question, then, is: can we prove an analogue of the
Andr\'asfai-Erd\H{o}s-S\'os theorem for general properties? Following the
terminology of Balogh, Bollob\'as and Simonovits~\cite{BalBolSim1}, who
investigated in detail the number and typical structure of $n$-vertex graphs in
general monotone properties, we make one further definition to facilitate
our quest.

\begin{definition}[Decomposition Family] Let $\mathcal{L}$ be a family of
non-empty graphs. Let $\mathcal{M}=\mathcal{M}(\mathcal{L})$ be the family of
minimal graphs $M$ with the property that, for some integer $t$, the graph $G$
obtained from $K_{p(\mathcal{L})}(t)$ by inserting a copy of $M$ into one part
is not $\mathcal{L}$-free. We call $\mathcal{M}$ the \emph{decomposition family}
of $\mathcal{L}$.
\end{definition}

One would hope that the r\^ole of the function $\biex$ in Theorem~\ref{AESext}
can be replaced by $\ex\big(n,\mathcal{M}(\mathcal{L})\big)$ for general
properties. Sadly this is not quite true: we can prove sharp results only when
$\mathcal{M}(\mathcal{L})$ is a finite set.

\begin{theorem} To any family $\mathcal{L}$ and finite subset $\mathcal{M'}$ of
$\mathcal{M}(\mathcal{L})$ there is associated a constant
$C=C(\mathcal{M'},\mathcal{L})$ such that whenever $\eps>0$ is given, there is
$n_0$ for which the following holds. Whenever $n\geq n_0$ and $G$ is an $n$-vertex graph with minimum degree exceeding
$\left(1-\frac{3}{3p(\mathcal{L})-1}+\eps\right)n$ which is $\mathcal{L}$-free,
then one can delete at most $C\ex(n,\mathcal{M'})$ edges from $G$
to obtain a $p(\mathcal{L})$-partite graph.
\end{theorem}

To see that, when $\mathcal{M}(\mathcal{L})=\mathcal{M}'$ is finite, this
theorem is best possible in the same sense as Theorem~\ref{AESext}, one need only repeat the two constructions at the end of
Section 1, this time with $\ex\big(n,\mathcal{M}(\mathcal{L})\big)$ replacing
$\biex(n,H)$. The proof of the theorem is essentially identical to the proof of
Theorem~\ref{AESext}: we reduce it to a brief sketch of the
required modification.

\begin{proof}
  Given $\mathcal{L}$ and a finite subset $\mathcal{M}'$ of
  $\mathcal{M}(\mathcal{L})$, we let $H$ be a graph in $\mathcal{L}$ with
  $\chi(H)=p(\mathcal{L})+1$ which amongst all such graphs minimises
  $\sigma(H)$. By definition and finiteness of $\mathcal{M}'$, there exists an
  integer $t$ such that for every $M\in\mathcal{M'}$, the graph obtained by inserting a copy of
  $M$ into one part of $K_{p(\mathcal{L})}(t)$ is not $\mathcal{L}$-free.
  
  We repeat the steps proving Theorem~\ref{AESext} for the selected graph
  $H\in\mathcal{L}$, with one alteration, at the point where
  Corollary~\ref{NiceCons} is used. At this point, we use instead the following
  statement.
  
  Given $\delta>0$, $p(\mathcal{L})$, $t$ and $\mathcal{M}'$, there exists a
  constant $C'$ with the following property. If $E$ is a set of edges in an
  $n$-vertex graph $G$ such that for every edge $uv$ in $E$, the edge $uv$ lies
  in at least $\delta n^{p(\mathcal{L})-1}$ copies of $K_{p(\mathcal{L})+1}$, then either $|E|\le C'\ex(n,\mathcal{M}')$, or $G$ contains a copy of
  $K_{p(\mathcal{L})}(t)$ in one of whose parts is a copy of some graph in
  $\mathcal{M}'$.
  
  The proof of this statement is a tiny modification to that
  of Corollary~\ref{NiceCons}. By the definition of $t$, if $G$ is an
  $\mathcal{L}$-free graph, then the second conclusion cannot hold, and thus we
  obtain $|E|\le C'\ex(n,\mathcal{M}')$, much as in the proof of
  Theorem~\ref{AESext}. It is at this step only that we
  require finiteness of $\mathcal{M}'$: the deduction that if $\ex(n,\mathcal{M}')$
  is smaller than $n-1$ then $\sigma(H)=1$ (which is the only other place in the
  proof of Theorem~\ref{AESext} that the function $\biex$ is used) does not
  require that $\mathcal{M}'$ be a finite set.
\end{proof}

One might think that there should be some way to avoid this finiteness
condition---but there is not. This is perhaps not so surprising: the task of
enumerating the $n$-vertex graphs in an $\mathcal{L}$-free property is closely
related to the extremal problem studied here, and it has already been shown by
Balogh, Bollob\'as and Simonovits~\cite{BalBolSim1} that when the decomposition
family is not finite, considering $\ex\big(n,\mathcal{M}(\mathcal{L})\big)$ may
not lead to sharp results.

Let $S_i$ be the graph obtained from the
cycle $C_{2i}$ and a disjoint complete bipartite graph $K_{i^4,i^4}$ by
inserting all edges between the $C_{2i}$ and $K_{i^4,i^4}$. Let
$\mathcal{S}=\{S_i:i\ge 2\}$. We claim that this family provides a
counterexample.

\begin{theorem} We have $\ex\big(n,\mathcal{M}(\mathcal{S})\big)\le 2n-2$, but
for every $C>0$, for all sufficiently large $n$, there is an $n$-vertex
$\mathcal{S}$-free graph with minimum degree at least
$\big(1-\tfrac{1}{p(\mathcal{S})}-\tfrac{1}{C}\big)n$ which cannot be made
$p(\mathcal{S})$-partite by removing $C(2n-2)$ edges.
\end{theorem}
\begin{proof}
Since every $S_i$ has chromatic number four, we have $p(\mathcal{S})=3$. It is
straightforward to verify that $\mathcal{M}(\mathcal{S})$ consists of all even cycles (and no other graphs).
Now suppose $F$ is an $\mathcal{M}(\mathcal{S})$-free $n$-vertex graph. Let $F'$
be a bipartite subgraph of $F$ with the maximum number of edges: then $F'$
cannot contain any cycles, and therefore $e(F')\le n-1$. Since $e(F)\le 2e(F')$,
we obtain $\ex\big(n,\mathcal{M}(\mathcal{S})\big)\le 2n-2$.

Given a constant $C$, we construct graphs as follows.

First, for all sufficiently large $n$, there exists a bipartite graph $G'$ on
$n/3$ vertices with $3Cn$ edges and girth at least $2\sqrt{\log n}$. To see
this, consider the random bipartite graph $\mathbf{H}$ with parts of size $n/6$
and edge probability $p=144C/n$. The expectation of $e(\mathbf{H})$ is
$4Cn$, and for each $g$, the expected number of cycles of length $g$ in
$\mathbf{H}$ is at most $(n/6)^g p^g=(24C)^g$. Let $X(\mathbf{H})$ count the
number of (even) cycles of length at most $2\sqrt{\log n}$ in $\mathbf{H}$. Then
by linearity of expectation we have
\[\mathbb{E}\big[e(\mathbf{H})-X(\mathbf{H})\big]\ge 4Cn-\sqrt{\log
n}(24C)^{2\sqrt{\log n}}\ge 3Cn\] for sufficiently large $n$. It follows that
there exists some $H$ with $e(H)-X(H)\ge 3Cn$; we can remove one edge from each
cycle of length at most $2\sqrt{\log n}$ to yield the desired graph $G'$ with at
least $3Cn$ edges and girth at least $2\sqrt{\log n}$.

Second, for all sufficiently large $n$, there exists a bipartite graph $G''$ on
$2n/3$ vertices with minimum degree $\big(1-\tfrac{1}{C}\big)n/3$ containing no copy of $K_{\log^2 n,\log^2
n}$. To see this, consider the random bipartite graph $\mathbf{H}'$ with parts
of size $n/3$ and edge probability $1-\tfrac{1}{2C}$. By a standard use of
Chernoff's inequality, a.a.s.\ every vertex of $\mathbf{H}'$ has degree at least
$\big(1-\tfrac{1}{C}\big)n/3$. Furthermore, the expected number of copies of
$K_{\log^2n,\log^2n}$ in $\mathbf{H}'$ is given by
\[ \binom{n/3}{\log^2 n}^2\big(1-\frac{1}{2C}\big)^{\log^4 n}<2^{2\log^3
n}\big(1-\frac{1}{2C}\big)^{\log^4 n}=o(1)\,, \]
and hence a.a.s.\ $\mathbf{H}'$ contains no copy of $K_{\log^2 n,\log^2 n}$: in
particular, our desired $G''$ exists for all sufficiently large $n$.

Finally, we let $G$ be obtained by taking the disjoint union of $G'$ and $G''$,
and inserting all edges between the two. We presume that $n>2^4$, so that $G'$ is $C_4$-free.

Observe that for each $i\ge 2$, the graph $S_i$ has chromatic number $4$, and
possesses exactly one $4$-partition. Furthermore, there is only one pair of
parts of $S_i$ whose induced bipartite graph may not contain $C_4$: namely the
pair inducing $C_{2i}$. Since $G$ is also $4$-partite, and the graph $G'$ which
makes up two of its partition classes is $C_4$-free, if we wish to find a copy
of $S_i$ in $G$ we have no choice but to start by finding a copy of $C_{2i}$ in
$G'$. This dooms us to failure when $i<\sqrt{\log n}$. On finding our copy
of $C_{2i}$, we are again left with no choice but to embed the remaining
$K_{i^4,i^4}$ of $S_i$ in $G''$: but this is impossible when $i\ge\sqrt{\log
n}$. It follows that for every $i$, $G$ is $S_i$-free, so $G$ is
$\mathcal{S}$-free. But we cannot make $G$ be $3$-partite without deleting at
least $e(G')\ge 3Cn>C(2n-2)$ edges, as required.
\end{proof}

This example is very similar to that given by Balogh, Bollob\'as and
Simonovits~\cite{BalBolSim1}. However, in their (tripartite) example, the even
cycles are joined completely to very large independent sets, and either a cycle
is too short to appear in some not too sparse graph (such as our $G'$, although
they mention several possible constructions) or the independent set is so huge
that $n$ vertices do not suffice to contain it. One feels that this is somehow
`cheating', and that it is preferable that the function $\ex(n,\mathcal{M})$
should be made small due to graphs in $\mathcal{L}$ which are on (preferably
much) less than $n$ vertices: hence the example given here. It is
straightforward to verify that excluding all even cycles of length up to $\log n$ from an $n$-vertex graph $G$
yields a linear upper bound on $e(G)$: in our example, the cycle $C_{2\log n}$
is contained in $\mathcal{M}(\mathcal{S})$ due to the graph $S_{\log n}$ on
$2\log n+2\log^4 n\ll n$ vertices; in the example of Balogh et al., the graph
responsible for excluding $C_{\log n}$ has more than $2^{\log^2 n}=n^{\log n}\gg
n$ vertices.

If the reader disagrees on this point, then the simpler example
in~\cite{BalBolSim1} works equally well here: if the reader too is disquieted at
the thought of using enormous graphs to make $\ex(n,\mathcal{M})$ small, then
perhaps it is comforting to observe that our example is also sufficient for
their conclusions: almost all bipartite graphs with parts of size $n/3$ do not
contain $K_{\log^2n,\log^2n}$.
\end{document}